\documentstyle[12pt]{article}

\newtheorem{theorem}{Theorem}
\newtheorem{lemma}{Lemma}

\newcommand{\be}{\begin{equation}}
\newcommand{\ee}{\end{equation}}

\begin{document}
\bibliographystyle{plain}

\thispagestyle{empty}
\setcounter{page}{0}

\vspace {2cm}

{\Huge Guszt\'av Morvai: 
 Guessing the Output of a Stationary Binary Time Series.}

\vspace {2cm}

{\Large In:  Foundations of statistical inference (Shoresh, 2000), pp. 207--215, Contrib. Statist., Physica, Heidelberg, 2003.
}

\vspace {2cm}

\
\begin{abstract}
The forward prediction problem for a binary   time series $\{X_n\}_{n=0}^{\infty}$ is to estimate the probability that $X_{n+1}=1$ based on the observations $X_i$, $0\le i\le n$ without prior knowledge of the distribution of the process $\{X_n\}$. 
It is known that this is not possible if one estimates at all values of $n$. 
 We present a simple procedure which will attempt to make such a prediction 
 infinitely often at carefully selected stopping times chosen 
by the algorithm. The growth rate of the stopping times is also exhibited.  
\end{abstract}

\section{Introduction}

T. Cover in \cite{Cover75} asked two fundamental questions concerning estimation for stationary and ergodic binary   processes. 
Cover's first question was as follows. 

\bigskip
\noindent
{\bf Question 1} {\it Is there an estimation scheme $f_{n+1}$ for the value $P(X_1=1|X_0,X_{-1},\dots,X_{-n})$ such that $f_{n+1}$ depends solely on the observed data segment $X_0,X_{-1},\dots, X_{-n}$ and 
$$
\lim_{n\to\infty}f_{n+1}(X_0,X_{-1},\dots,X_{-n})-P(X_{1}=1|X_0,X_{-1},\dots,X_{-n})=0
$$
almost surely for all stationary and ergodic binary   time series $\{X_n\}$?. 
}

\bigskip
\noindent
This question was answered by Ornstein \cite{Ornstein78} by constructing such a scheme. (See also Bailey \cite{Bailey76}.)
Ornstein's scheme is not a simple one and  the proof of consistency is rather sophisticated. A much simpler scheme and proof of consistency were provided by Morvai, Yakowitz, Gy\"orfi \cite{MoYaGy96}. (See also Weiss \cite{Weiss00}.)

\noindent
Here is Cover's second question. 

\bigskip
\noindent
{\bf Question 2} {\it Is there an estimation scheme $f_{n+1}$ for the value $P(X_{n+1}=1|X_0,X_1,\dots,X_n)$ such that $f_{n+1}$ depends solely on the data segment $X_0,X_1,\dots, X_n$ and 
$$
\lim_{n\to\infty}f_{n+1}(X_0,X_1,\dots,X_n)-P(X_{n+1}=1|X_0,X_1,\dots,X_n)=0
$$
almost surely for all stationary and ergodic binary    time series $\{X_n\}$?. }

\bigskip
\noindent
This question was answered by Bailey \cite{Bailey76} in a negative way, that is, he showed that there is no such scheme. 
(Also see Ryabko \cite{Ryabko88}, Gy\"orfi,  Morvai, Yakowitz \cite{GYMY98}
 and  Weiss \cite{Weiss00}.) Bailey used the technique of cutting and stacking developed by Ornstein \cite{Ornstein74} (see also Shields \cite{Shields91}). Ryabko's construction was  based on a function of an infinite state Markov-chain. 
This negative result can be interpreted as follows. Consider a weather forecaster  whose task it is to predict the probability of the event 'there will be rain tomorrow' given the observations up to the present day. Bailey's result says that the difference between the estimate and the true conditional probability cannot   eventually be small for all stationary weather processes. The difference will be big infinitely often. 
These results show that there is a great  difference between Questions 1 and 2.  Question 1 was addressed by Morvai, Yakowitz, Algoet \cite{MoYaAl97} and a very simple estimation scheme was given which satisfies the statement in Question 1 in probability instead of almost surely. 
Now consider a less ambitious goal than  Question 2: 

\bigskip
\noindent
{\bf Question 3} {\it 
Is there a sequence of stopping times $\{\lambda_n\}$ and an estimation scheme 
$f_n$ which  depends on the observed data segment $(X_0,X_1,\dots,X_{\lambda_n})$ such that 
$$
\lim_{n\to\infty} \left( f_n(X_0,X_1,\dots,X_{\lambda_n})-
P(X_{\lambda_n+1}=1|X_0,X_1,\dots,X_{\lambda_n})\right) =0
$$
almost surely for all stationary binary time series $\{X_n\}$?
}

\bigskip
\noindent
It turns out that the  answer is affirmative and such a scheme will be exhibited below.
This result can be interpreted as if the weather forecaster can refrain from predicting, that is, he may say that he does not want to predict today, but will predict at infinitely many time instances, and the difference between the prediction and the true conditional probability will vanish  almost surely at the stopping times.


\section{Forward Estimation for Stationary Binary   Time Series}

Let $\{X_n\}_{n=-\infty}^{\infty}$ denote a two-sided stationary binary   time series. For $n\ge m$, it will be convenient to use the notation 
$X^n_m=(X_m,\dots,X_n)$.   
For $k=1,2,\ldots$,  define the sequences $\{\tau _k\}$ and  $\{\lambda_k\}$ 
recursively. Set $\lambda_0=0$. 
Let
$$
{\tau}_k=
\min\{t>0 : X_{t}^{\lambda_{k-1}+t}=X_0^{\lambda_{k-1}}\}
$$
and
$$
\lambda_k=\tau_k+\lambda_{k-1}.
$$
(By stationarity,   
the string $X_{0}^{\lambda_{k-1}}$ must appear  in the sequence
$X_1^{\infty}$ almost surely. )
The $k$th estimate of $P(X_{\lambda_k+1}=1|X_0^{\lambda_k})$ is denoted  by $P_k,$
and is defined as 
\begin{equation}
\label{fkdistrestimate1}
P_k=
{1\over k-1}\sum_{j=1}^{k-1} X_{\lambda_j+1} 
\end{equation} 
For an arbitrary stationary binary   time series $\{Y_n\}_{n=-\infty}^0$,  
for  $k=1,2,\ldots$,  define the sequence ${\hat \tau} _k$ and ${\hat\lambda_k}$
recursively. Set ${\hat\lambda}_0=0$. 
Let
$$
{\hat\tau}_k=
\min\{t>0 : {Y}_{-\hat\lambda_{k-1}-t}^{-t}={Y}_{{-\hat\lambda}_{k-1}}^{0}\}
$$
and let
$$
{\hat\lambda}_k={\hat\tau}_k+{\hat\lambda}_{k-1}.
$$
When there is ambiguity as  to which time series $\hat\tau_k$ and $\hat\lambda_k$ are to be applied, we will use the notation
$\hat\tau_k(Y^0_{-\infty})$ and $\hat\lambda_k(Y^0_{-\infty})$.

It will be useful to define another time series $\{{\tilde X}_n\}_{n=-\infty}^0$ as
\begin{equation}\label{deftildeprocess}
{\tilde X}_{-\lambda_k}^{0}:=X_0^{\lambda_k} \ \ \mbox{for all $k\ge 1$.}
\end{equation}
Since $X_{\lambda_{k+1}-\lambda_k}^{\lambda_{k+1}}=X_{0}^{\lambda_{k}}$ the above definition is correct.
Notice that it is immediate that $\hat\tau_k({\tilde X}_{-\infty}^0)=\tau_k$ and  $\hat\lambda_k({\tilde X}_{-\infty}^0)=\lambda_k$.

\begin{lemma} \label{distrlemma}
The two time series $\{{\tilde X}_n\}_{n=-\infty}^0$ and $\{X_n\}_{n=-\infty}^{\infty}$ have identical  distribution, that is, for all $n\ge 0$, and $x^0_{-n}\in\{0,1\}^{n+1}$,
$$
P({\tilde X}^0_{-n}=x^0_{-n})=P(X^0_{-n}=x^0_{-n}).
$$ 
\end{lemma}
{\sc Proof}
First we prove that    
\begin{equation}\label{lambdaeq}
P({\tilde X}_{-n}^0=x_{-n}^0,{\hat\lambda}_k({\tilde X}_{-\infty}^0)=n)=
P(X_{-n}^0=x_{-n}^0,{\hat\lambda}_k(X_{-\infty}^0)=n).
\end{equation}
Indeed, by~(\ref{deftildeprocess}),  
${\tilde X}_{-\hat\lambda_k({\tilde X}_{-\infty}^0)}^0=X_0^{\lambda_k}$, and it yields  
$$
P({\tilde X}_{-n}^0=x_{-n}^0,{\hat\lambda}_k({\tilde X}_{-\infty}^0)=n)=
P(X_0^n=x_{-n}^0,\lambda_k=n),
$$
and by stationarity, 
$$
P(X_0^n=x_{-n}^0,\lambda_k=n)=
P(X_{-n}^0=x_{-n}^0,\hat\lambda_k(X_{-\infty}^0)=n)
$$
and (\ref{lambdaeq}) is proved.
Apply ~(\ref{lambdaeq})  in order to get 
\begin{eqnarray*}
\lefteqn{P({\tilde X}^0_{-n}=x^0_{-n}) }\\
&=& \sum_{j=n}^{\infty} P({\tilde X}^0_{-n}=x^0_{-n},
{\hat\lambda}_n({\tilde X}_{-\infty}^0)=j) \\
&=& \sum_{j=n}^{\infty} \sum_{x^{-n-1}_{-j} \in\{0,1\}^{j-n}} P({\tilde X}^0_{-j}=x^0_{-j}, {\hat\lambda}_n({\tilde X}_{-\infty}^0)=j) \\
&=& \sum_{j=n}^{\infty} \sum_{x^{-n-1}_{-j}\in\{0,1\}^{j-n}} P(X^0_{-j}=x^0_{-j}, {\hat\lambda}_n(X_{-\infty}^0)=j) \\
&=& \sum_{j=n}^{\infty} P(X^0_{-n}=x^0_{-n},
{\hat\lambda}_n(X_{-\infty}^0)=j) \\
&=& P(X^0_{-n}=x^0_{-n})
\end{eqnarray*}
and Lemma~\ref{distrlemma} is proved. 

\bigskip
\noindent
Since $\{X_n\}_{n=-\infty}^{\infty}$ is a stationary time series, by Lemma~\ref{distrlemma} so is $\{ {\tilde X}_n\}_{n=-\infty}^0$. Since a stationary time series can  always be extended to be a two-sided time series we have also defined   
$\{{\tilde X}_n\}_{n=-\infty}^{\infty}$.  Now we prove the universal consistency of the estimator $P_k$. 

\begin{theorem} \label{Theorem1} For all stationary binary time series $\{X_n\}$ and estimator defined in 
(\ref{fkdistrestimate1}),
\begin{equation}
\label{fkstatement}
\lim_{k\to\infty} \left( P_k- P(X_{\lambda_k+1}=1|X_0^{\lambda_k})\right) =0\ \ \mbox{almost surely.}
\end{equation}
Moreover,
\begin{equation}
\label{limit}
\lim_{k\to\infty} P_k= \lim_{k\to\infty} P(X_{\lambda_k+1}=1|X_0^{\lambda_k})=P({\tilde X}_1=1|{\tilde X}^0_{-\infty})
\end{equation}
almost surely.
\end{theorem}
{\sc Proof}
\begin{eqnarray*}
\lefteqn{P_k -P(X_{\lambda_k+1}=1|X_0^{\lambda_k})}\\
&=& {1\over k-1}\sum_{j=1}^{k-1}
 \{X_{\lambda_j+1}-P(X_{\lambda_j+1}=1|X_0^{\lambda_j})]\}\\
&+& {1\over k-1}\sum_{j=1}^{k-1}
 \{P(X_{\lambda_j+1}=1|X_0^{\lambda_j})
-P(X_{\lambda_k+1}=1|X_0^{\lambda_k})\}\\
&=& {1\over k-1}\sum_{j=1}^{k-1} \Gamma_j+
{1\over k-1}\sum_{j=1}^{k-1} (\Delta_j-\Delta_k).
\end{eqnarray*}
Observe that
 $\{\Gamma_j,\sigma(X_0^{\lambda_j+1})\}$ is a bounded martingale difference sequence 
for $1\leq j<\infty$. To see this note that
$\sigma(X_0^{\lambda_j+1})$ is monotone
increasing, and  $\Gamma_j$   is measurable with respect to 
$\sigma(X_0^{\lambda_j+1})$, and $E(\Gamma_j|X_0^{\lambda_{j-1}+1})=0$ for $1\leq j<\infty$.
Now apply  Azuma's exponential
bound for bounded martingale differences in Azuma \cite{Azuma67} to
get that for any $\epsilon>0$,
$$P\left(\left|{1\over (k-1)}\sum_{j=1}^{k-1} \Gamma_j\right|>\epsilon\right)
\le 2\exp(-\epsilon^2 (k-1)/2).$$
After summing the right side over $k$, and appealing to
the Borel-Cantelli lemma for a sequence of $\epsilon$'s tending to zero we get  
$${1\over (k-1)}\sum_{j=1}^{k-1} \Gamma_j \to 0\ \ \mbox{almost surely.}$$
It remains to show 
$${1\over k-1}\sum_{j=1}^{k-1} \Delta_j-\Delta_k\to 0\ \ \mbox{almost surely.}$$
Define 
$$
p_{k,n}(x_{-n}^0)=P(X_{\lambda_k+1}=1|X_0^{\lambda_k}=x_{-n}^0,\lambda_k=n)
$$
and (applying $\hat\lambda_k$ to the time series $\{ {\tilde X}_n\}_{n=-\infty}^0$) 
$$
{\tilde p}_{k,n}(x_{-n}^0)=P({\tilde X}_1=1|{\tilde X}_{-\hat\lambda_k}^0=x_{-n}^0,\hat\lambda_k=n).
$$
Now  the fact that $\lambda_k=\hat\lambda_k$ and  Lemma~\ref{distrlemma} together imply  
\begin{equation} \label{equalprob1}
p_{k,n}(x_{-n}^0)={\tilde p}_{k,n}(x_{-n}^0).
\end{equation}
By~ (\ref{deftildeprocess}) and (\ref{equalprob1}),
\begin{equation}\label{equalprob2}
{\tilde p}_{k,\lambda_k}(X_{\lambda_k}^0)=
{\tilde p}_{k,\hat\lambda_k}({\tilde X}_{-\hat\lambda_k}^0).
\end{equation}
Combine (\ref{equalprob1}) and (\ref{equalprob2}) in order  
to get  
$$
P(X_{\lambda_k+1}=1|X_0^{\lambda_k})=
P({\tilde X}_1=1|{\tilde X}_{-\hat\lambda_k}^0).
$$
Notice that $\{P({\tilde X}_1=1|{\tilde X}_{-\hat\lambda_k}^0),\sigma({\tilde X}_{-\hat\lambda_k}^0)\}$ is a bounded martingale and so it converges almost surely to $P({\tilde X}_1=1|{\tilde X}^0_{-\infty})$, and so does 
$P(X_{\lambda_k+1}=1|X_0^{\lambda_k})$. We have proved that $\Delta_j$ converges almost surely. Now Toeplitz lemma yields that  
$
{1\over k-1}\sum_{j=1}^{k-1} (\Delta_j-\Delta_k) \to 0
$ almost surely.  
The proof of Theorem~\ref{Theorem1} is complete.

\bigskip 

\section{The Growth Rate of the Stopping Times}

\bigskip
\noindent
The next result shows that the growth of  the stopping times $\{\lambda_k\}$ is rather rapid. 
Let $p(x^0_{-n})=P(X^0_{-n}=x^0_{-n})$. 

\begin{theorem} \label{FiniteAsLemmaComplexity} 
Let $\{X_n\}$ be a stationary and ergodic binary time series. Suppose 
that $H>0$ where 
$$H=\lim_{n\to\infty}-{1\over n+1} E \log p(X_0,\dots,X_n)$$
is the process entropy. 
Let $0<\epsilon<H$ be arbitrary. 
Then for $k$ large enough,  
\begin{equation}
\label{complexityAs}
{\lambda}_{k}(\omega)\ge c^{ c^{ {\cdot}^{ {\cdot}^c}}} \ 
\mbox{almost surely,}
\end{equation}
where the height of the tower is $k-K$, $K(\omega)$ is a finite number 
which depends on $\omega$, and $c=2^{H-\epsilon}$. 
\end{theorem}
{\sc Proof}
Since by (\ref{deftildeprocess}), $\lambda_k={\hat\lambda}_k({\tilde X}_{-\infty}^0)$, and by Lemma~\ref{distrlemma} the time series 
$\{X_n\}_{-\infty}^{\infty}$ and $\{ {\tilde X}_n\}_{-\infty}^{\infty}$ have identical  distributions, and hence the same entropy, it is enough to prove the result for $\hat\lambda_k({\tilde X}_{-\infty}^0)$. Now $\hat\tau_k$ and $\hat\lambda_k$ are evaluated on the process $\{\tilde X_n\}_{n=-\infty}^0$.  
\noindent
For $0<l<\infty$ define 
$$
R(l)=\min\{j\ge l+1: {\tilde X}_{-l-j}^{-j}={\tilde X}_{-l}^{0} \}.
$$
By Ornstein and Weiss \cite{OrWe93}, 
\begin{equation}
\label{OrnWeissentr}
{1\over l+1} \log R(l) \to H \ \mbox{almost surely.}
\end{equation}

\noindent
First we show that if $H>0$ then for $k$ large enough 
${\hat\tau}_{k+1}> {\hat\lambda}_k$ almost surely. 
We argue by contradiction. Suppose  that ${\hat\tau}_{k+1}\to \infty$ and 
${\hat\tau}_{k+1}\le {\hat\lambda}_k$ infinitely often. 
Then 
$$
{\tilde X}_{-\hat\lambda_k}^{0}={\tilde X}_{-\hat\lambda_k-\hat\tau_{k+1}}^{-\hat\tau_{k+1}}
$$
and $\hat\tau_{k+1}\le \hat\lambda_k$ infinitely often. Hence 
$$
{\tilde X}_{-\hat\tau_{k+1}+1}^{0}={\tilde X}_{-\hat\tau_{k+1}-\hat\tau_{k+1}+1}^{-\hat\tau_{k+1}} 
$$
infinitely often 
and $R(\hat\tau_{k+1}-1)\le \hat\tau_{k+1}$ infinitely often. 
Then by~(\ref{OrnWeissentr}), 
\begin{eqnarray*}
H &=&\lim_{k\to\infty} {1\over \hat\tau_{k+1}}\log R(\hat\tau_{k+1}-1)\\
&\le& \lim_{k\to\infty} {1\over \hat\tau_{k+1}}\log \hat\tau_{k+1}\\
&=& 0
\end{eqnarray*}
provided that $\hat\tau_k\to\infty$. 
Now assume that $\eta=\sup_{0<k<\infty} \hat\tau_k$ is finite. 
Then $R(n\eta-1)=n\eta$. Now by~(\ref{OrnWeissentr}), 
\begin{eqnarray*}
 H &=&\lim_{n\to\infty} {1\over n\eta}\log R(n\eta-1)\\
&\le& \lim_{n\to\infty} {1\over n\eta}\log(n\eta)\\
&=& 0. 
\end{eqnarray*}
We have shown  that $H>0$ implies that for $k$ large enough 
$\hat\tau_{k+1}> \hat\lambda_k$ almost surely and hence 
for $k$ large enough $R(\hat\lambda_k)=\hat\tau_{k+1}$ almost surely. 
Hence by~(\ref{OrnWeissentr}), 
$$
{1\over \hat\lambda_k+1}\log\hat\tau_{k+1}\to H \ \mbox{almost surely.}
$$
Thus for almost every $\omega\in\Omega$ 
there exists a positive finite integer $K(\omega)$ such that for 
$k\ge K(\omega)$,  
${1\over \hat\lambda_k+1}\log\hat\tau_{k+1}>H-\epsilon$ 
and 
$$
\hat\lambda_{k+1}>\hat\tau_{k+1}>c^{\hat\lambda_k} \ \mbox{for $k\ge K(\omega)$}
$$
and the proof of Theorem~\ref{FiniteAsLemmaComplexity} is complete. 

\bigskip

\section{Guessing the Output at Stopping Time Instances}

\bigskip
If the weather forecaster is pressed to say simply will it rain or not tomorrow then we need a guessing scheme, rather than a predictor. 
Define the guessing scheme $\{{\bar X}_{\lambda_k}\}$ for the values $\{X_{\lambda_k+1}\}$ as 
$$
{\bar X}_{\lambda_k}=1_{\{P_k\ge 0.5\}}.
$$
Let $X^*_{\lambda_k}$ denote the Bayes rule, that is, 
$$
X^*_{\lambda_k}=1_{\{P(X_{\lambda_k+1}=1|X_0^{\lambda_k}) \ge 0.5\}}.
$$

\begin{theorem} \label{guessintheorem} Let $\{X_n\}_{n=-\infty}^{\infty}$ be a stationary binary  time series. 
The proposed guessing scheme ${\bar X}_{\lambda_k}$ works in the average at stopping times $\lambda_k$ just as well as the Bayes rule, that is, 
\begin{equation}   
\lim_{n\to\infty} \left( {1\over n} \sum_{k=1}^n 1_{\{\bar X_{\lambda_k}= X_{\lambda_k+1}\}}-
{1\over n} \sum_{k=1}^n 1_{\{X^*_{\lambda_k}= X_{\lambda_k+1}\}}\right) =0
\end{equation}
almost surely. 
Moreover, 
\begin{equation}
\lim_{k\to\infty} \left( P({\bar X}_{\lambda_k}= X_{\lambda_k+1}|X_0^{\lambda_k})-
 P(X^*_{\lambda_k}= X_{\lambda_k+1}|X_0^{\lambda_k})\right)
=0
\end{equation}
almost surely. 
\end{theorem}
{\sc Proof}

\begin{eqnarray*}
\lefteqn{ \sum_{k=1}^n 1_{\{{\bar X}_{\lambda_k}= X_{\lambda_k+1}\}}-
{1\over n} \sum_{k=1}^n 1_{\{X^*_{\lambda_k}= X_{\lambda_k+1}\}}=}\\
&& {1\over n} \sum_{k=1}^n \left[ 1_{\{ {\bar X}_{\lambda_k}= X_{\lambda_k+1}\}}-
P({\bar X}_{\lambda_k}= X_{\lambda_k+1}|X_0^{\lambda_k})\right]\\
&-&
{1\over n} \sum_{k=1}^n \left[1_{ \{X^*_{\lambda_k} = X_{\lambda_k+1} \}}-
P(X^*_{\lambda_k}= X_{\lambda_k+1}|X_0^{\lambda_k})\right]\\
&+&
{1\over n} \sum_{k=1}^n
\left[ P({\bar X}_{\lambda_k}= X_{\lambda_k+1}|X_0^{\lambda_k})-
P(X^*_{\lambda_k}= X_{\lambda_k+1}|X_0^{\lambda_k})\right] \\
&=& \Gamma_n+\Theta_n+\Psi_n.
\end{eqnarray*}
Now
$\Gamma_n$ and $\Theta_n$ tend to zero since they are averages of bounded martingale differences (cf. Azuma \cite{Azuma67}). Concerning the third term $\Psi_n$, it is enough to prove that 
$$
\lim_{k\to\infty} \left( P({\bar X}_{\lambda_k}= X_{\lambda_k+1}|X_0^{\lambda_k})-
P(X^*_{\lambda_k}= X_{\lambda_k+1}|X_0^{\lambda_k})\right) =0
$$
almost surely. 
To see this recall  the result in  Theorem~\ref{Theorem1},
$$\lim_{k\to\infty} P_k=\lim_{k\to\infty} 
P(X_{\lambda_k+1}=1|X_0^{\lambda_k})=P({\tilde X}_1=1|{\tilde X}^0_{-\infty})
$$
almost surely, and apply this in order to get 
\begin{eqnarray*}
\lefteqn{
\lim_{k\to\infty} [P({\bar X}_{\lambda_k}= X_{\lambda_k+1}|X_0^{\lambda_k})
-
P(X^*_{\lambda_k}= X_{\lambda_k+1}|X_0^{\lambda_k})]=}\\
&\lim_{k\to\infty}& \{ [P(P({\tilde X}_1=1|{\tilde X}^0_{-\infty})\neq 0.5, {\bar X}_{\lambda_k}= X_{\lambda_k+1}|X_0^{\lambda_k})\\
&-&
P(P({\tilde X}_1=1|{\tilde X}^0_{-\infty})\neq 0.5, X^*_{\lambda_k}= X_{\lambda_k+1}|X_0^{\lambda_k})]\\
&+&
[P(P({\tilde X}_1=1|{\tilde X}^0_{-\infty})= 0.5, {\bar X}_{\lambda_k}= X_{\lambda_k+1}|X_0^{\lambda_k})\\
&-&
P(P({\tilde X}_1=1|{\tilde X}^0_{-\infty})= 0.5, X^*_{\lambda_k}= X_{\lambda_k+1}|X_0^{\lambda_k})] \} \\
&=& 0.
\end{eqnarray*}
The proof of  Theorem~\ref{guessintheorem} is now complete.

\bigskip
\noindent   
{\bf Acknowledgments.}
The author wishes to thank Benjamin Weiss for helpful discussions and suggestions. 
This paper has been written by the auspices of the Hungarian National E\"otv\"os Fund. 
(Ez a cikk a Magyar \'Allami E\"otv\"os \"Oszt\"ond\'\i j t\'amogat\'as\'aval k\'esz\"ult.)

\end{document}